\documentclass[10pt]{amsart}
\usepackage[cp1251]{inputenc}
\usepackage[english,russian]{babel}
\usepackage{amsmath}
\usepackage{amssymb}
\usepackage{amsfonts}
\usepackage{graphicx}

\theoremstyle{definition}

\theoremstyle{plain}

\begin{document}
\renewcommand{\refname}{References}

\thispagestyle{empty}

\title{Resolvability and complete accumulation points}
\author{{A.E. Lipin}}%
\address{Anton Evgenevich Lipin
\newline\hphantom{iii} N.N. Krasovskii Institute of Mathematics and Mechanics
of the Ural Branch of the Russian Academy of Sciences,
\newline\hphantom{iii} Sofia Kovalevskaya street, 16,
\newline\hphantom{iii} 620108, Ekaterinburg, Russia}%
\email{tony.lipin@yandex.ru}%

\thanks{{\sc Lipin, A.E.,
Resolvability and complete accumulation points.} \\
The work was performed as part of research conducted in the Ural Mathematical
Center with the financial support of the Ministry of Science and Higher Education of
the Russian Federation.}
\thanks{\copyright \ 2023 A.E. Lipin}

\maketitle 
{\small
\begin{quote}
\noindent{\sc Abstract. } We prove that:

I. For every regular Lindel{\"o}f space $X$ if $|X|=\Delta(X)$ and $\mathrm{cf}|X|\ne\omega$, then $X$ is maximally resolvable;

II. For every regular countably compact space $X$ if $|X|=\Delta(X)$ and $\mathrm{cf}|X|=\omega$, then $X$ is maximally resolvable.

Here $\Delta(X)$, the dispersion character of $X$, is the minimum cardinality of a nonempty open subset of $X$.

Statements I and II are corollaries of the main result: for every regular space $X$ if $|X|=\Delta(X)$ and every set $A\subseteq X$ of cardinality $\mathrm{cf}|X|$ has a complete accumulation point, then $X$ is maximally resolvable.

Moreover, regularity here can be weakened to $\pi$-regularity, and the Lindel{\"o}f property can be weakened to the linear Lindel{\"o}f property.

\medskip

 \end{quote}
}

\section{Introduction}

Let us start with recalling the main definition (introduced by E.~Hewitt in 1943 \cite{Hewitt} and expanded by J.G.~Ceder in 1964 \cite{Ceder}).

\vskip 10pt

\noindent\textbf{1.1. Definition.}
Let $\kappa$ be a cardinal. A topological space $X$ is called $\kappa$-{\it resolvable} if $X$ contains $\kappa$ pairwise disjoint dense subsets.
The space $X$ is called {\it maximally resolvable} if $X$ is $\Delta(X)$-resolvable, where
$$\Delta(X)=\min\{|U| : U \text{ is nonempty open in } X\}$$
is called the dispersion character of $X$. The space $X$ is called {\it resolvable} if it is $2$-resolvable and {\it irresolvable} if it is crowded (i.e. without isolated points) but not resolvable.

\vskip 10pt

In the general case the question of $\kappa$-resolvability of some space $X$ is nontrivial for $\kappa$ from $2$ to $\Delta(X)$ (in particular, only for crowded spaces). Hewitt explored resolvability and irresolvability, the parameter $\kappa$ was added by Ceder. Hewitt found a method to construct irresolvable spaces and proved that all spaces of certain "nice"{} classes (metric spaces, $T_2$ compact spaces and some others) are resolvable \cite{Hewitt}. In 1964 J.G.~Ceder showed that actually all spaces in these classes are maximally resolvable \cite{Ceder}.

We refer the reader to a selective survey \cite{survey} made by O.~Pavlov in 2007 for more detailed information on resolvability. Some topological properties allow spaces to be irresolvable and some others imply resolvability but not maximal resolvability. Moreover, there are examples of $n$-resolvable spaces which are not $(n+1)$-resolvable \cite{Feng} and even spaces which are $\lambda$-resolvable for all $\lambda<\kappa$ but not $\kappa$-resolvable, where the cardinal $\kappa$ is uncountable regular \cite{J-S-Sz-01} (but it is impossible if $\mathrm{cf}(\kappa)=\omega$ \cite{Illanes}, \cite{Rao}).

Generalizations of compactness usually entail resolvability, although in many cases we still do not know if they entail maximal resolvability. One of the first studied generalizations is a Hausdorff compact generated space (i.e. a space, whose topology is defined by its compact subspaces). In 1975 N.V.~Velichko proved that every such space is resolvable \cite{Velichko}, in 1979 E.G.~Pytkeev proved that it is maximally resolvable \cite{Pytkeev}.

Resolvability of a Lindel{\"o}f space has been studied since at least 1998. The natural restriction here is to consider only spaces with uncountable dispersion character, because all countable spaces are Lindel{\"o}f but not all countable spaces are resolvable \cite{Hewitt}. In 1998 V.I.~Malykhin constructed an example of irresolvable Hausdorff Lindel{\"o}f space with uncountable dispersion character and asked whether a space $X$ is resolvable if $X$ is regular Lindel{\"o}f and $\Delta(X)\geq\omega_1$ \cite{Malykhin}.

In 2002 Pavlov proved that if $\Delta(X)\geq\omega_2$ then the space $X$ is $\omega$-resolvable \cite{Pavlov}. The case $\Delta(X)=\omega_1$ was open for two more years, until in 2004 M.A.~Filatova proved that $X$ is resolvable \cite{Filatova}. This answered Malykhin's question, but raised a natural question of $3$-resolvability. In 2007 and 2012  I.~Juhasz, L.~Soukup and Z.~Szentmiklossy have strengthened results of Pavlov and Filatova and finally proved that $X$ is $\omega$-resolvable \cite{J-S-Sz-0}, \cite{J-S-Sz-1}.  Also they proved that if, in addition, $|X|=\Delta(X)=\omega_1$ then $X$ is even $\omega_1$-resolvable. In \cite{J-S-Sz-1} they ask if such space $X$ is maximally resolvable.

Besides that, in recent years the same results were received for Lindel{\"o}f generated spaces, i.e. spaces with topology defined by Lindel{\"o}f subspaces (Filatova and A.V.~Osipov \cite{F-O}, Juhasz, Soukup and Szentmiklossy \cite{J-S-Sz-3}).

Another class of spaces whose maximal resolvability remains generally unknown is the class of regular countably compact spaces (note that Malykhin constructed an example of an irresolvable Hausdorff countably compact too \cite{Malykhin}). E.G.~Pytkeev proved that every such a space is $\omega_1$-resolvable \cite{Pytkeev-2}, then in 2016-2017 some results were received for more wide class of pseudocompact spaces (van Mill \cite{Mill}, Ortiz-Castillo and Tomita \cite{O-T}, Juhasz, Soukup and Szentmiklossy \cite{J-S-Sz-2}; at least the last one gives even continuum-resolvability but only for spaces with small enough cellularity).

Lindel{\"o}f and countably compact spaces have strong properties connected with complete accumulation points.

Recall that a point $x$ of a space $X$ is called a {\it complete accumulation} point for an infinite set $A\subseteq X$ if for every neighborhood $U$ of $x$ we have $|U\cap A|=|A|$. We denote $A^\circ$ the set of all complete accumulation points of the set $A$ (as in \cite{J-S-Sz-1}). Following the terminology of \cite{J-S-Sz-1}, we say that a space $X$ is $\kappa$-compact if every set $A\in[X]^\kappa$ has a complete accumulation point. 

It is well-known (and easy to prove) that all countably compact spaces are $\omega$-compact and all Lindel{\"o}f spaces are $\kappa$-compact for every cardinal $\kappa$ with uncountable cofinality (actually this is a criterion for wider class of linearly Lindel{\"o}f spaces). In this paper we investigate resolvability of regular $\kappa$-compact spaces to receive corollaries (5.1 and 5.2) for countably compact and (linearly) Lindel{\"o}f spaces. The main result of this paper is Theorem 4.6.

We also note that actually throughout the article regularity can be weakened to $\pi$-regularity (which means that for every nonempty open set $U$ there is a nonempty open set $V$ such that $\overline{V}\subseteq U$).

\section{Preliminaries}

We assume the following notation and conventions.

\begin{itemize}

\item If we say that $A=\{x_\alpha : \alpha<\kappa\}$ is {\it enumeration} we mean that for $\alpha\ne\beta$ always $x_\alpha\ne x_\beta$.

\item And if we say that $A=\{x_\alpha : \alpha<\kappa\}$ is {\it indexing} then we do not mean such a condition.

\item Symbol $\bigsqcup$ denote {\it disjoint union} in the following sense: its equal to the usual union, but using it we assume that the united sets are pairwise disjoint.

\item We use letters $\kappa,\lambda$ for cardinals and $\alpha,\beta,\gamma,\delta$ and in one case $\xi$ for ordinals.

\item $\mathrm{cf}(\kappa)$ is the cofinality of the cardinal $\kappa$.

\item If $S$ is a set and $\kappa$ is a cardinal then $[S]^\kappa$ (resp. $[S]^{<\kappa}$) is the family of all sets $A\subseteq S$ such that $|A|=\kappa$ ($|A|<\kappa$).

\item {\it Space} means topological space.

\item $\Delta(X)$ is the dispersion character of the space $X$, i.e. minimum cardinality of a nonempty open subset of $X$.

\item If $A$ is a subset of some space then $A^\circ$ is the set of all complete accumulation points of $A$.


\item A space $X$ is $\kappa$-{\it compact} if for all $A\in[X]^\kappa$ we have $A^\circ\ne\emptyset$.

\item If $A$ is a subset of some space then $\overline{A}$ is the closure of $A$.

\end{itemize}

\section{Set-theoretical Propositions}

Three propositions of this section must be well-known, but we give proofs for the sake of accuracy.

\vskip 10pt

\noindent\textbf{3.1. Definition.}
Let $\kappa$ be an infinite cardinal and $S$ be a set. A family $\mathcal{E}$ of subsets of $S$ is said to be $\kappa$-{\it almost disjoint} on $S$ if for all $A\in \mathcal{E}$ we have $|A|\geq\kappa$ and for all different $A,B\in\mathcal{E}$ we have $|A\cap B|<\kappa$.

\vskip 10pt

\noindent\textbf{3.2. Proposition.}
\textit{ Let $S$ be a set of regular cardinality $\kappa$ and a family $\mathcal{E}\subseteq[S]^{\kappa}$ be $\kappa$-almost disjoint and such that $|\mathcal{E}|\leq\kappa$. Then it is possible to choose sets $A^*\subseteq A$ for all $A\in\mathcal{E}$ in such a way that $|A\setminus A^*|<\kappa$ and the family $\{A^* : A\in\mathcal{E}\}$ is disjoint.}

\vskip 5pt

\noindent\textit{Proof.}
Let us choose an arbitrary enumeration $\mathcal{E}=\{A_\alpha : \alpha<\lambda\}$ where $\lambda\leq\kappa$.
Now for all $\alpha<\lambda$ we define $A_\alpha^* := A_\alpha \setminus \bigcup\limits_{\beta<\alpha}A_\beta$. It is clear that all conditions are satisfied.

\begin{flushright}$\square$\end{flushright}

\vskip 10pt

\noindent\textbf{3.3. Proposition.}
\textit{ Let $S$ be a set of cardinality $\kappa\geq\omega$ and a family $\mathcal{E}\subseteq[S]^{<\kappa}$ satisfy the condition $|\mathcal{E}|\leq\kappa$. Then there is a set $B\in[S]^{\mathrm{cf}(\kappa)}$ such that for any $A\in S$ we have $|A\cap B|<\mathrm{cf}(\kappa)$.}

\vskip 5pt

\noindent\textit{Proof.}
The case of regular $\kappa$ is trivial. Assume that $\mathrm{cf}(\kappa)<\kappa$.

We claim that there is an indexing $\mathcal{E}=\{A_\alpha : \alpha<\kappa\}$ with the property that for each $\gamma<\kappa$ always
$|\bigcup\limits_{\alpha<\gamma}A_\alpha|<\kappa$ (recall that by using the word "indexing"{} we mean that one element can receive more than one index).

Let us show that such an indexing exists. At first we take arbitrary indexing $\{B_\alpha : \alpha<\kappa\}=\mathcal{E}$ such that every element of $\mathcal{E}$ has $\kappa$ indices. We also choose any $C\in\mathcal{E}$. Now we define $A_\alpha$ to be $B_\alpha$ if $|B_\alpha|\leq|\alpha|$ and $C$ otherwise. Obviously, $|\bigcup\limits_{\alpha<\gamma}A_\alpha|\leq|\gamma|+|C|<\kappa$.

Now we take any increasing $\mathrm{cf}(\kappa)$-sequence $(\gamma_\xi : \xi<\mathrm{cf}(\kappa))$ with supremum $\kappa$. Finally, for all $\xi<\mathrm{cf}(\kappa)$ we choose arbitrary $x_\xi\in S\setminus\bigcup\limits_{\alpha<\gamma_\xi}A_\alpha$.
It is clear that the set $B := \{x_\xi : \xi<\mathrm{cf}(\kappa)\}$ satisfies all conditions.

\begin{flushright}$\square$\end{flushright}

\vskip 10pt

\noindent\textbf{3.4. Proposition.}
\textit{ If $S$ is a set of cardinality $\kappa\geq\omega$ then there is a $\mathrm{cf}(\kappa)$-almost disjoint family $\mathcal{E}$ on $S$ such that $|\mathcal{E}|>\kappa$.}

\vskip 5pt

\noindent\textit{Proof.}
Let us take any maximal $\mathrm{cf}(\kappa)$-almost disjoint family $\mathcal{E}$ with $\kappa$ or more elements and prove that actually $|\mathcal{E}|>\kappa$.
Assume, on the contrary, that $|\mathcal{E}|=\kappa$.

If $\kappa$ is singular, then by Proposition 3.3 we can append to $\mathcal{E}$ some new set of cardinality $\mathrm{cf}(\kappa)$ keeping $\kappa$-almost disjointness. This is a contradiction.

If $\kappa$ is regular, then by Proposition 3.2 we can choose subsets $A^*\subseteq A$ for all $A\in\mathcal{E}$ in such a way that $|A\setminus A^*|<\kappa$ and the family $\{A^* : A\in\mathcal{E}\}$ is disjoint.
Let us take any elements $x_A \in A^*$ for all $A\in\mathcal{E}$ and construct $B := \{x_A : A\in\mathcal{E}\}$.
Clearly, $|B|=\kappa$ and we can append $B$ to the family $\mathcal{E}$ keeping $\kappa$-almost disjointness. This is a contradiction again.

\begin{flushright}$\square$\end{flushright}

\section{Main result}

Our goal here is to prove Theorem 4.6. Our proof splits into two cases, singular and regular. Now we are ready to prove the singular case. In terms of \cite{Pavlov} it can be reduced to the fact that, under the conditions of the theorem, $(\kappa,\mathrm{cf}(\kappa))$-trace is dense in the space $X$ (see Definition 1.2 and Proposition 2.1 in \cite{Pavlov}), but we give a full version of the proof.

\vskip 10pt

\noindent\textbf{4.1. Statement.}
\textit{ Let $X$ be a regular space such that $|X|=\Delta(X)=\kappa>\mathrm{cf}(\kappa)$ and $X$ is $\mathrm{cf}(\kappa)$-compact. Then $X$ is maximally resolvable.}

\vskip 5pt

\noindent\textit{Proof.}
Let us denote $M$ the set of all points $x\in X$ such that there is $R_x\in[X]^{<\kappa}$ with the property that for every $A\in[X\setminus R_x]^{\mathrm{cf}(\kappa)}$ we have $x\notin A^\circ$. We also denote $H := X\setminus M$.

First of all, we prove that $H$ is dense. Let $U$ be any nonempty open set in $X$. We denote $\mathcal{E} := \{R_x\cap\overline{U} : x\in M\}$. By Proposition 3.3 applied to $\mathcal{E}$ and $\overline{U}$ there is a set $B\in[\overline{U}]^{\mathrm{cf}(\kappa)}$ such that for all $x\in M$ we have $|R_x\cap\overline{U}\cap B|<\mathrm{cf}(\kappa)$ and hence $x\notin B^\circ$ (otherwise we would have $x\in (B\setminus R_x)^\circ$ what is a contradiction with the definition of $M$). So nonempty set $B^\circ$ is contained in $H\cap\overline{U}$. Considering that the space $X$ is regular, this means that $H$ is dense.

Now we choose any indexing $H=\{x_\alpha : \alpha<\kappa\}$ with the property that every point $x\in H$ has $\kappa$ indices. Recursively by $\alpha<\kappa$ then by $\beta<\alpha$ we choose sets $A_{\alpha,\beta}\in[X]^{\mathrm{cf}(\kappa)}$ in such a way that $x_\alpha\in A_{\alpha,\beta}^\circ$ and if $(\alpha,\beta)\ne(\gamma,\delta)$ then $A_{\alpha,\beta}\cap A_{\gamma,\delta}=\emptyset$.

For all $\beta<\kappa$ we define $D_\beta := \bigsqcup\limits_{\beta<\alpha<\kappa}A_{\alpha,\beta}$. It is easy to see that sets $D_\beta$ are pairwise disjoint and for all $x\in H$ we have $x\in\overline{D_\beta}$, so $\overline{D_\beta}\supseteq H$, hence $\overline{D_\beta}=X$.

\begin{flushright}$\square$\end{flushright}

The regular case of Theorem 4.6 requires a little more work. We prove it with the help from the following notion (compare our technique with Definition 2.5 and Theorem 2.7 in \cite{J-S-Sz-1}. We mainly repeat their idea with little technical complication which allows us to use Proposition 3.4 and prove the cornerstone Lemma 4.3).

\vskip 10pt

\noindent\textbf{4.2. Definition.}
Let $\kappa$ be a cardinal and $X$ be a space. We say that a collection $\Sigma$ of families of subsets of the space $X$ is a $\kappa$-{\it fission} of the space $X$ at a nonempty set $H \subseteq X$ if all of the following conditions are met:
\begin{enumerate}

\item[(A)] $|\Sigma|=\kappa$, and if $\mathcal{S}\in\Sigma$, then $|\mathcal{S}|\leq\kappa$, and if $A\in\mathcal{S}\in\Sigma$, then $|A|=\kappa$;

\vskip 10pt

\item[(B)] if $A\in\mathcal{S}\in\Sigma$, $B\in\mathcal{T}\in\Sigma$ and $(\mathcal{S},A)\ne(\mathcal{T},B)$, then $|A\cap B|<\kappa$;


\vskip 10pt

\item[(C)] if $\mathcal{S}\in\Sigma$, then $\bigcup\limits_{A\in\mathcal{S}}A^\circ=H$.

\end{enumerate}

\noindent We say that $X$ is $\kappa$-{\it fissile} at a set $H\subseteq X$ if there is a $\kappa$-fission of $X$ at $H$.

\vskip 10pt

Let us prove that such an object really exists.

\vskip 10pt

\noindent\textbf{4.3. Lemma.}
\textit{ Let $X$ be a $\kappa$-compact space, where $\kappa=|X|$ is regular. Then there is a nonempty set $H\subseteq X$ such that $X$ is $\kappa$-fissile at $H$.}

\vskip 5pt

\noindent\textit{Proof.}
By Proposition 3.4 there is some $\kappa$-almost disjoint family $\mathcal{E}$ on $X$ such that $|\mathcal{E}|>\kappa$.
For every $x\in X$ we denote $\mathcal{L}_x$ the family $\{A\in\mathcal{E} : x\in A^\circ\}$.
Let $M$ be the set of all points $x\in X$ such that $|\mathcal{L}_x|\leq\kappa$.
Define $H := X \setminus M$ and $\mathcal{F} := \mathcal{E} \setminus \bigcup\limits_{x\in M}\mathcal{L}_x$.

Let us prove that $H$ is not empty.
Since $|M|\leq|X|=\kappa$ and for all $x\in M$ we have $|\mathcal{L}_x|\leq\kappa$, it follows that $|\mathcal{F}|=|\mathcal{E}|>\kappa$.
Choose any $A\in\mathcal{F}$. For all $x\in M$ we have $A\notin\mathcal{L}_x$, which by definition means $x\notin A^\circ$. Consequently, the nonempty set $A^\circ$ is contained in $H$. So $H\ne\emptyset$.

Note that for every $x\in H$ we have more than $\kappa$ sets $A\in\mathcal{F}$ such that $x\in A^\circ$ (actually $\kappa$ would be enough).
It allows us to choose for all $\alpha<\kappa$ and $x\in H$ pairwise different sets $A_{x,\alpha}\in\mathcal{F}$ such that $x\in A_{x,\alpha}^\circ$.
For every $\alpha<\kappa$ we define $\mathcal{S}_\alpha := \{A_{x,\alpha} : x\in H\}$.
Finally, denote $\Sigma := \{\mathcal{S}_\alpha : \alpha<\kappa\}$. It is easy to check that $\Sigma$ is a $\kappa$-fission of $X$ at $H$.

\begin{flushright}$\square$\end{flushright}

\vskip 10pt

\noindent\textbf{4.3.1. Corollary.}
\textit{ Let $X$ be a regular $\kappa$-compact space, where $\kappa=|X|=\Delta(X)$ is regular. Then there is a disjoint family $\mathcal{H}$ of subsets of $X$ such that $\bigsqcup\mathcal{H}$ is dense and the space $X$ is $\kappa$-fissile at every $H\in\mathcal{H}$.}

\vskip 5pt

\noindent\textit{Proof.}
Let us denote $\mathcal{H}$ an arbitrary maximal disjoint family of subsets of $X$ such that the space $X$ is $\kappa$-fissile at every $H\in\mathcal{H}$.
We prove that $\mathcal{H}$ is the required family.

The only thing left to prove is that $\bigsqcup\mathcal{H}$ is dense.
Let $U$ be a nonempty open set in $X$.
By Lemma 4.3 the subspace $\overline{U}$ contains some subset $K$ such that there is a $\kappa$-fission $\Sigma$ of $\overline{U}$ at $K$.
Easy to see that $\Sigma$ is also $\kappa$-fission of $X$ at $K$,
hence by maximality of $\mathcal{H}$ there is $H\in\mathcal{H}$ such that $H\cap K\ne\emptyset$.
Consequently, $\bigsqcup\mathcal{H}\cap\overline{U}\ne\emptyset$.
Since $X$ is regular, we can deduce now that the set $\bigsqcup\mathcal{H}$ is dense.

\begin{flushright}$\square$\end{flushright}

\vskip 10pt

\noindent\textbf{4.4. Lemma.}
\textit{ Let $X$ be a regular $\kappa$-compact space, where $\kappa=|X|=\Delta(X)$ is regular. Then there is a dense set $H\subseteq X$ such that the space $X$ is $\kappa$-fissile at $H$.}

\vskip 5pt

\noindent\textit{Proof.}
We choose any family $\mathcal{H}$ with the properties from Corollary 4.3.1. Let us show that $\bigsqcup\mathcal{H}$ is the required set.
For every $H\in\mathcal{H}$ we denote $\Sigma_H$ any $\kappa$-fission of $X$ at $H$.

We claim that for any $H,K\in\mathcal{H}$, $A\in\mathcal{S}\in\Sigma_H$ and $B\in\mathcal{T}\in\Sigma_K$
if $(H,\mathcal{S},A)\ne(K,\mathcal{T},B)$, then $|A\cap B|<\kappa$.
If $H=K$, then it follows from property (B) of definition of a $\kappa$-fission.
Let $H\ne K$.
Assume, on the contrary, that $|A\cap B|=\kappa$. It follows that $(A\cap B)^\circ\ne\emptyset$ and also $(A\cap B)^\circ\subseteq H\cap K$, hence we have a contradiction with the condition of disjointness of $\mathcal{H}$.

Let us choose arbitrary enumerations $\Sigma_H=\{\mathcal{S}_{H,\alpha} : \alpha<\kappa\}$. 
For all $\alpha<\kappa$ we denote $\mathcal{S}_\alpha$ the family $\bigsqcup\limits_{H\in\mathcal{H}}\mathcal{S}_{H,\alpha}$.
Finally, we construct $\Sigma := \{\mathcal{S}_\alpha : \alpha<\kappa\}$.
It is easy to check that the collection $\Sigma$ is a $\kappa$-fission of $X$ at $\bigsqcup\mathcal{H}$.

\begin{flushright}$\square$\end{flushright}

\vskip 10pt

\noindent\textbf{4.5. Lemma.}
\textit{ Let a space $X$ be $\kappa$-fissile at a dense set $H\subseteq X$, where $\kappa$ is regular. Then the space $X$ is $\kappa$-resolvable.}

\vskip 5pt

\noindent\textit{Proof.}
Let a collection $\Sigma$ be a $\kappa$-fission of the space $X$ at the set $H$.
By Proposition 3.2 we can choose pairwise disjoint sets $A^*\subseteq A$ for all $A\in\mathcal{S}\in\Sigma$ in such a way that $|A\setminus A^*|<\kappa$. Note that $(A^*)^\circ=A^\circ$.

For every $\mathcal{S}\in\Sigma$ we denote $D_\mathcal{S} := \bigsqcup\limits_{A\in\mathcal{S}}A^*$. The sets $D_\mathcal{S}$ are pairwise disjoint and for all $\mathcal{S}\in\Sigma$ we have $\overline{D_\mathcal{S}}\supseteq\overline{H}=X$. Since $|\Sigma|=\kappa$, it follows that $X$ is $\kappa$-resolvable.

\begin{flushright}$\square$\end{flushright}

Connecting Lemmas 4.4 and 4.5 together, we complete the proof of the regular case of our theorem. We are ready to formulate our main result in the general form.

\vskip 10pt

\noindent\textbf{4.6. Theorem.} 
\textit{ Let $X$ be a regular $\mathrm{cf}(\kappa)$-compact space, where $\kappa=|X|=\Delta(X)$. Then $X$ is maximally resolvable.}

\vskip 10pt

Let us end this section by noting once again that the condition of regularity can be weakened to $\pi$-regularity in all statements here (and in the next section too).

\vskip 10pt

\section{Corollaries}

It is well-known that if a space $X$ is countably compact, then every infinite subset $A\subseteq X$ has a limit point, hence $X$ is $\omega$-compact. Consequently, we have

\vskip 10pt

\noindent\textbf{5.1. Corollary.}
\textit{ Let $X$ be a regular countably compact space, $|X|=\Delta(X)$ and $\mathrm{cf}|X|=\omega$. Then $X$ is maximally resolvable.}

\vskip 10pt

We recall that a space $X$ is said to be {\it linearly Lindel{\"o}f} if every linearly ordered by inclusion open cover of $X$ contains a countable subcover. It is known that the space $X$ is linearly Lindel{\"o}f if and only if $X$ is $\lambda$-compact for every infinite $\lambda$ with uncountable cofinality. Every Lindel{\"o}f space is obviously linearly Lindel{\"o}f. We formulate the next corollaries for Lindel{\"o}f spaces only, but note that the same is true for all linearly Lindel{\"o}f spaces.

\vskip 10pt

\noindent\textbf{5.2. Corollary.}
\textit{ Let $X$ be a regular Lindel{\"o}f space, $|X|=\Delta(X)$ and $\mathrm{cf}|X|\ne\omega$. Then $X$ is maximally resolvable.}

\vskip 10pt

Let us finish with a little remark on the condition $|X|=\Delta(X)$. In \cite{Elkin} A.G.~Elkin observed the following useful fact:

\vskip 5pt

\noindent $\bullet$ {\it A space $X$ is $\kappa$-resolvable if and only if every nonempty open set in $X$ contains a nonempty $\kappa$-resolvable subspace.}

\vskip 5pt

In a regular Lindel{\"o}f space $X$ every nonempty open set contains a nonempty Lindel{\"o}f subspace $Y$ such that $\Delta(X)\leq\Delta(Y)=|Y|\leq|X|$. Therefore, investigation of (maximal) resolvability of regular Lindel{\"o}f spaces is reduced to the case $|X|=\Delta(X)$. Clearly, the same is true for countable compactness and any other property inherited by closed subspaces. In particular, Corollary 5.2 has its own

\vskip 10pt

\noindent\textbf{5.2.1. Corollary.}
\textit{ Let $X$ be a regular Lindel{\"o}f space. Suppose that for every cardinal $\kappa$ such that $\Delta(X)\leq\kappa\leq|X|$ it follows that $\mathrm{cf}(\kappa)\ne\omega$ (for instance, this is true if $\omega<\Delta(X)\leq|X|<\aleph_\omega$). Then $X$ is maximally resolvable.}

\vskip 10pt

Now the question if a regular Lindel{\"o}f space $X$ with an uncountable dispersion character is maximally resolvable is reduced to the case $|X|=\Delta(X)>\mathrm{cf}|X|=\omega$. This case remains unsolved, actually we do not even know if such space is uncountably resolvable.

\section{Acknowledgements}

The author is grateful to Maria~A. Filatova for constant attention to this work and to Vladislav~R. Smolin and Vladimir~V. Ivchenko for their help with editing.

\end{document}